\documentstyle[amsfonts,12pt]{article}
\newcommand{\lon}{\longrightarrow}
\newcommand{\rar}{\rightarrow}

\newcommand{\rH}{\mbox{H}}
\newcommand{\Z}{{\Bbb Z}}

\newcommand{\C}{{\Bbb C}}
\newcommand{\R}{{\Bbb R}}
\newcommand{\ot}{\otimes}
\newcommand{\Id}{\mbox{Id}}
\newcommand{\Beq}{\begin{equation}}
\newcommand{\Eeq}{\end{equation}}
\newcommand{\Beqr}{\begin{eqnarray}}
\newcommand{\Eeqr}{\end{eqnarray}}
\newcommand{\Beqrn}{\begin{eqnarray*}}
\newcommand{\Eeqrn}{\end{eqnarray*}}
\newcommand{\Ba}{\begin{array}}
\newcommand{\Ea}{\end{array}}
\newcommand{\Barr}{\begin{array}}
\newcommand{\Earr}{\end{array}}
\newcommand{\Bi}{\begin{itemize}}
\newcommand{\Ei}{\end{itemize}}
\newcommand{\Bc}{\begin{center}}
\newcommand{\Ec}{\end{center}}
\newcommand{\fg}{{\frak g}}
\newcommand{\fh}{{\frak h}}
\newcommand{\f}{{\cal O}}

\newcommand{\cC}{{\cal C}}
\newcommand{\cM}{{\cal M}}
\newcommand{\cA}{{\cal A}}
\newcommand{\cB}{{\cal B}}

\newcommand{\cK}{{\cal K}}

\newcommand{\al}{\alpha}

\newcommand{\ga}{\gamma}
\newcommand{\var}{\varepsilon}
\newcommand{\la}{\lambda}

\newcommand{\tv}{\tilde{v}}
\newcommand{\ta}{\tilde{a}}
\newcommand{\tl}{\tilde{l}}

\newcommand{\tk}{\tilde{k}}

\newcommand{\tj}{\tilde{j}}

\newcommand{\vst}{\vspace{2 mm}}

\newcommand{\vse}{\vspace{8 mm}}
\begin{document}

\title{$L_{\infty}$-algebra of an unobstructed\\
deformation functor}

\author{ S.A.\ Merkulov}
\date{}
\maketitle

\begin{abstract}
 This is a comment on the Kuranishi method of constructing
analytic deformation spaces. It is based on a simple observation
that the Kuranishi map can always be inverted in the category of
$L_{\infty}$-algebras. The $L_{\infty}$-structure obtained by
this inversion is used to define an ''unobstructed'' deformation
functor which is always representable by a smooth pointed moduli
space.

The singular nature of the original Kuranishi deformation
space emerges in this setting merely as a result of the truncation of
this ``naive'' $L_{\infty}$-algebra controlling the deformations
to a usual differential Lie algebra.
\end{abstract}

\section{Introduction}
\paragraph{1.0.} Typical moduli problems in geometry are plagued with
obstructions and associated singularities. It is noticed in this paper
that, in the situations when the Kuranishi method of
constructing an analytic deformation space applies, there is
always a ``naive'' way of bypassing (rather than overcoming)
the obstructions by suitably extending the deformation problem
from the categories of associative or Lie algebras to the
associated strong homotopy versions. The resulting moduli spaces are
always smooth.

\vst

The singular nature of the original Kuranishi deformation
space emerges in this setting merely as a result of the truncation of the
``non-obstructed" \footnote{more precisely, quasi-isomorphic to a differential
Abelian Lie algebra} $L_{\infty}$-algebra controlling the deformations
to a usual differential Lie algebra.

\vst

The paper is based on  a simple observation that, roughly speaking, the Kuranishi
map can
always be inverted in the category of $L_{\infty}$-algebras.

\vst

To make the above statements more clear, we shall remind in the next subsection a few
facts about the Deformation Theory, and then formulate the main results again.

\paragraph{1.1. The deformation functor.} The standard approach
to constructing the analytic deformation space of a given mathematical
structure $\cA$\, consists of two key steps \cite{D,Ko,Ku1,Ku2,GM1,GM2}:
\Bi
\item[1)] Associate to $\cA$ a ``controlling" differential graded
Lie algebra $(\fg=\bigoplus_{k\in \Z} \fg^k, d, [\, ,\, ])$ over a field $k$
(which is usually $\R$ or $\C$) and define the deformation functor
$$
\Ba{rccc}
\mbox{Def}^{0}_{\fg}: & \left\{\Ba{l} \mbox{the category of Artin}\\
                    \mbox{$k$-local algebras $\cB$}\\
                    \mbox{with maximal ideals $m_{\cB}$}\Ea \right\}&
\lon & \left\{\mbox{the category of sets}\right\}
\Ea
$$
as follows
$$
\mbox{Def}^{0}_{\fg}(\cB)=\left\{\Gamma\in C(\fg\ot m_{\cB})^1 \mid d\Gamma
+ \frac{1}{2}[\Gamma, \Gamma]=0\right\},
$$
where $C(\fg\ot m_{\cB})^1$ is a complement to the 1-coboundaries,
$d(\fg \ot m_{\cB})^0$, in $(\fg\ot m_{\cB})^1$ and $\cB$ is viewed as a
$\Z$-graded
algebra concentrated in degree zero (so that $(\fg\ot m_{\cB})^i=\fg^i\ot m_{\cB}$).
\item[2)] Obtain a Hodge theory on $\fg$ and apply the Kuranishi
construction to represent the deformation functor by a germ, $\f_{p}$,
of the structure sheaf on a pointed analytic space $(\cK^0_{\fg}, p\in \cK^0_{\fg})$.
\Ei

The basic examples are the differential graded algebras $(TM\ot
\Omega^{0,\bullet}M, \bar{\partial}, \mbox{Schouten bracket})$
and $(E\ot E^* \ot \Omega^{0,\bullet}M, \bar{\partial}, \mbox{standard
bracket})$. The first controls deformations of a given complex structure on a
smooth manifold $M$, while the second controls deformations
of holomorphic structures on a
given complex  vector bundle $E\rar M$.

\vst

The tangent space, $\mbox{Def}_{\fg}^{0}(k[\var])$, to the functor
$\mbox{Def}_{\fg}^{0}$ is isomorphic to the first cohomology group $\rH^1(\fg)$
of the complex $(\fg,d)$. If one extends in the obvious way the
above deformation functor
to the category of arbitrary $\Z$-graded $k$-local Artin algebras (which may not
be concentrated in degree 0), one gets the functor
$\mbox{Def}_{\fg}^{\Z}$
with the tangent space isomorphic to  the full cohomology group
$\rH^*(\fg)$. Moreover, if $(\fg, d, [\, , \, ])$ happens to
be formal, one can  construct an associated Kuranishi moduli
space $\cK^{\Z}_{\fg}$ representing  $\mbox{Def}_{\fg}^{\Z}$
\cite{GM2}.

\vst

This extended deformation functor $\mbox{Def}_{\fg}^{\Z}$
has been used recently, in the mirror
symmetry context, by Barannikov
and Kontsevich \cite{BaKo,Ba} in constructing the smooth extended moduli
space of complex structures on a Calabi-Yau manifold. According to
Kontsevich's homological mirror symmetry conjecture \cite{Ko2}, the ``extended
holomorphic deformations" of a complex vector bundle on a
Calabi-Yau manifold may also play an important role in mirror symmetry
(cf.\ \cite{Vafa}). The extended  moduli space of special Lagrangian
submanifolds in a Calabi-Yau manifold was constructed in \cite{Me}.

%%%%%%%%%%%%%%%%%%%%%%%%%%%%%%%%%%%%%%%%%%%%%%%%%%%%%%%%%%
\paragraph{1.2. Obstructions.} If the Lie algebra
$(\fg,d, [\, , \, ])$ is such that $\rH^2(\fg)\neq 0$,
 then the associated deformation functor $\mbox{Def}^{0}_{\fg}$ is usually
obstructed in the sense that the local moduli space $\cK^0_{\fg}$ has
singularities.
\vst

On the other hand, if the Lie algebra
$(\fg,d, [\, , \, ])$ is  formal but the induced
Lie bracket on its cohomology,
\Beq
[\, , \, ]: \rH^*(\fg) \times \rH^*(\fg) \lon \rH^*(\fg), \label{obs}
\Eeq
is non-vanishing, the same problem plagues the extended deformation functor
$\mbox{Def}^{\Z}_{\fg}$ --- the associated Kuranishi moduli space
$\cK^{\Z}_{\fg}$ is  singular.

%%%%%%%%%%%%%%%%%%%%%%%%%%%%%%%%%%%%%%%%%%%%%%%%%%%%%%%%%%%
\paragraph{1.3. Bypassing the obstructions.}
It is shown in Sect.\ 3 of this paper
that the Hodge theory on  $(\fg, d, [\, , \, ])$ gives rise
canonically to the structure of a $L_{\infty}$-algebra on the vector
space $\fg$, i.e.\ to a set of linear maps $\mu_n: \Lambda^n \fg
\rar \fg[2-n]$ satisfying the higher Jacobi
identities. We use this structure to define a  modified deformation functor,
$$
\mbox{MDef}_{\fg}^{0}(\cB) = \left\{\Gamma \in C(\fg\ot
m_{\cB})^1 \mid \sum_{k=1}^{\infty} \frac{(-1)^{k(k+1)/2}}{k!}
\mu_k(\Gamma, \ldots, \Gamma)=0\right\},
$$
which associates to a concentrated in degree zero $\Z$-graded
Artin\footnote{our apologies for this  awkward series of
adjectives, but its truncation leads to a very different object,
see below.} algebra $\cB$
with the maximal ideal
$m_{\cB}$ a set of solutions to the Maurer-Cartan equations in the
$L_{\infty}$ algebra $(\fg\ot m_{\cB}, \mu_*)$. Here $C(\fg\ot m_{\cB})^1$
is a  complement to $d((\fg\ot m_{\cB})^0)$ in
$(\fg\ot m_{\cB})^1$ (alternatively, we may replace above the subspace
$C(\fg\ot m_{\cB})^1$ by
$(\fg\ot m_{\cB})^1$ and take the quotient
by a natural gauge equivalence; this change is possible because $(\fg, \mu_*)$
is quasi-isomorphic to a differential Abelian Lie algebra).

\vst

Analogously, one defines the extended deformation functor
$\mbox{MDef}_{\fg}^{\Z}$ (by simply omitting the words
{\em concentrated in degree zero}\, in the sentence above).

\vst

The main results of the paper are
\Bi
\item[(i)] The functor $\mbox{MDef}_{\fg}^{0}$
is  {\em always}\, representable by an
analytic germ, $\f_{p}$, of a {\em smooth}\, pointed moduli space $({\cM}^0_{\fg},
p\in {\cM}^0_{\fg})$;
the tangent space to ${\cM}^0_{\fg}$ at $p$ is canonically isomorphic to
$\rH^1(\fg)$.

\item [(ii)] The original Kuranishi moduli space $\cK^0_{\fg}$ space arises
in this setting as an analytic
subspace in ${\cM}^0_{\fg}$ where both functors, the classical
one,
$\mbox{Def}_{\fg}^{0}$,   and the
modified one, $\mbox{MDef}_{\fg}^{0}$ ,  agree\footnote{With $\mbox{MDef}^0$ one can
canonically associate a truncated deformation functor, $\mbox{mDef}^0(\cB):=
\left\{\Gamma\in C(\fg\ot m_{\cB})^1 \mid d\Gamma
+ \frac{1}{2}\mu_2(\Gamma, \Gamma)=0, \ \mu_{n\geq 3}(\Gamma, \ldots, \Gamma)=0\right\}
$, which is representable by a pointed analytic subspace $(M^0_{\fg},p)\subset
(\cM_{\fg}^0,p)$; the word {\em agree}\, in the text means that
$\mbox{mDef}_{\fg}^{0}\simeq
\mbox{Def}_{\fg}^{0}$ and $M^0_{\fg}\simeq \cK^0_{\fg}.$}.

If the obstructions
vanish, $\rH^2(\fg)=0$, then $\mbox{MDef}_{\fg}^{0}\simeq
\mbox{Def}_{\fg}^{0}$ and $\cK^0_{\fg}\simeq \cM^0_{\fg}$.

\item[(iii)] The functor $\mbox{MDef}_{\fg}^{\Z}$
is  {\em always}\, representable by an
analytic germ, $\f_{p}$, of a pointed {\em smooth}\,moduli superspace $({\cM}^{\Z}_{\fg},
p\in {\cM}^{\Z}_{\fg})$;
the tangent space to ${\cM}^{\Z}_{\fg}$ at $p$ is canonically isomorphic to
$\rH^*(\fg)$.

\item[(iv)] If the differential Lie superalgebra $(\fg, d, [\, , \, ])$
is formal and the obstruction map (\ref{obs}) vanishes,
then $\mbox{MDef}_{\fg}^{\Z} \simeq \mbox{Def}_{\fg}^{\Z}$
and ${\cM}^{\Z}_{\fg}\simeq \cK^{\Z}_{\fg}$.

\Ei

The paper is organized as follows. In Sect.\ 2 we give a brief
introduction into the theory of $L_{\infty}$-algebras and
list a few necessary facts
about Maurer-Cartan equations. In Sect.\ 3 we
construct a particular $L_{\infty}$-algebra used in the
definition of the functor
and prove the main results (i)-(iv) stated above.

\vst

A few words about notations. The category $\Z$-graded
vector spaces over a field $k$ contains an object $[1]=\bigoplus_{i\in
\Z}[1]^i$ defined by
$$
[1]^{i}=\left\{\Ba{ll} k & \mbox{if}\ i=-1,\\
0 & \mbox{if}\ i\neq-1.\Ea\right.
$$
Its tensor powers are denoted by $[n]$, and the tensor product of a graded
vector space $\fg=\bigoplus_{i\in \Z} \fg^i$ with $[n]$ is denoted by
$\fg[n]$. A homogeneous vector $v\in \fg$ viewed as an element of
$\fg[n]$  is denoted by $v[n]$.

\vst

For a homogeneous element $v\in \fg^i\subset \fg$ we write $\tilde{v}:=
i\bmod 2\Z \in \Z_2$. Analogously, for an integer $n\in \Z$ we write $\tilde{n}=n
\bmod 2\Z$.

\vst

If one is interested in the extended deformations functors
$\mbox{Def}_{\fg}^{\Z}$ or $\mbox{MDef}_{\fg}^{\Z}$ only, the $\Z$-grading
above and below can be safely replaced by the associated $\Z_2$-grading --- no
essential information will be lost.

%%%%%%%%%%%%%%%%%%%%%%%%%%%%%%%%%%%%%%%%%%%%%%%%%%%%%%%%
%%%%%%%%%%%%%%%%%%%%%%%%%%%%%%%%%%%%%%%%%%%%%%%%%%%%%%%%
\section{ Strong homotopy algebras and Maurer-Cartan equations}

\paragraph{2.1. $L_{\infty}$-algebras.}
A {\em strong homotopy Lie algebra}, or shortly $L_{\infty}$-{\em algebra},
is by definition a $\Z$-graded vector space $\fg$ equipped with linear maps,
$$
\Ba{rcccc}
\mu_k: & \Lambda^k \fg & \lon & \fg[2-k] & \\
& v_1\ot\ldots \ot v_k & \lon & \mu_k(v_1, \ldots, v_k), & \ \ \ \ k\geq 1,
\Ea
$$
satisfying, for any $n\geq 1$ and
 any $v_1, \ldots, v_n \in V$,
 the following
{\em higher order Jacobi identities},
\Beq \label{id}
\sum_{k+l=n+1}\sum_{\sigma\in Sh(k,n)} (-1)^{\tilde{\sigma}+k(l-1)} e(\sigma;v_1,\ldots,v_n)
 \mu_l\left(\mu_k(v_{\sigma(1)},\ldots,v_{\sigma(k)}),
v_{\sigma(k+l)},\ldots, v_{\sigma(n)}\right)=0,
\Eeq
where $Sh(k,n)$ is the set of all permutations $\sigma:\{1, \ldots, n\} \rar
\{1,\ldots,n\}$ which satisfy $\sigma(1)<\ldots< \sigma(k)$ and $\sigma(k+1)<\ldots
<\sigma(n)$. The symbol
$e(\sigma; v_1, \ldots, v_n)$ (which we abbreviate from now on to $e(\sigma)$)
 stands for the {\em Koszul
sign}\,  defined by the equality
$$
v_{\sigma(1)}\wedge \ldots v_{\sigma(n)}= (-1)^{\tilde{\sigma}}e(\sigma)
v_1\wedge \ldots \wedge v_n,
$$
$\tilde{\sigma}$ being the parity of the permutation
$\sigma$.

\vst

This notion and the associated notion of  $A_{\infty}$-algebra (reminded
below) are due to Stasheff \cite{St1,St2,St3}.

\vst

The first three higher order Jacobi identities have the form
\begin{description}
\item[$n=1$:]\ \  $d^2=0$,\\
\item[$n=2$:]\ \  $d[v_1, v_2]= [dv_1, v_2] + (-1)^{\tv_1} [v_1,
dv_2]$,\\
\item[$n=3$:]\ \ $[[v_1, v_2],v_3] +
(-1)^{(\tv_1+\tv_2)\tv_3}[[v_3, v_1],v_2]
+ (-1)^{\tv_1(\tv_2+\tv_3)}[[v_2, v_3],v_1]$ =
$-d\mu_3(v_1,v_2,v_3) - \mu_3(dv_1,v_2,v_3) - (-1)^{\tv_1}\mu_3(v_1,dv_2,v_3)
-(-1)^{\tv_1+\tv_2}\mu_3(v_1,v_2,dv_3)$,\\
\end{description}
where we denoted $dv_1:=\mu_1(v_1)$ and $[v_1, v_2] := \mu_2(v_1,v_2)$.

\vst

Therefore $L_{\infty}$-algebras with $\mu_k=0$ for $k\geq 3$
are nothing but the usual graded differential Lie algebras with the
differential $\mu_1$ and the Lie bracket given by
$\mu_2$. If, furthermore, $\mu_1=0$, one gets  the class of usual graded Lie
algebras.

\vst

%%%%%%%%%%%%%%%%%%%%%%%%%%%%%%%%%%%%%%%%%%%%%%%%%%%%%%%%%%%

\paragraph{2.2. Another (conceptually better) definition of $L_{\infty}$ algebra.}
Let $W$ be a $\Z$-graded vector space and let $\odot^*W=\bigoplus_{n=1}^{\infty}\odot^n W$
be the associated symmetric tensor algebra with its induced $\Z$-grading.
We equip $\odot^*W$ with the structure of cosymmetric {\em
coalgebra}\, by setting
$$
\Delta(w_1\odot\ldots \odot w_n) =\sum_{i=1}^n \sum_{\sigma\in
Sh(i,n)} e(\sigma) \left(w_{\sigma(1)}\odot \ldots \odot
w_{\sigma(i)}\right)
\otimes \left(w_{\sigma(i+1)}\odot \ldots \odot
w_{\sigma(n)}\right).
$$

\noindent\mbox{{\bf 2.2.1. Fact}~\cite{LS}}. A $L_{\infty}$-algebra
structure on a graded vector space $\fg$ is equivalent to a
codifferential on the coalgebra $(\odot^*(\fg[1]), \Delta)$, i.e.\ to a
linear map
$$
Q: \odot^* (\fg[1]) \lon \left(\odot^* (\fg[1])\right)[1]
$$
satisfying the conditions
\begin{itemize}
\item[(i)] $\Delta \circ Q=\left(Q\otimes \Id + \Id\otimes Q\right)\circ
\Delta$,\\
\item[(ii)] $Q^2=0$.
\end{itemize}
The first condition simply says that $Q$ is a coderivation of
the coalgebra $(\odot^* (\fg[1]), \Delta)$ and hence is completely
determined by the compositions (``values in cogenerators"),
$$
\hat{\mu}_n:\odot^n (\fg[1]) \lon \odot^* (\fg[1]) \stackrel{Q}{\lon}
\left(\odot^* (\fg[1])\right)[1]
\lon \fg[2],
$$
for all $n\geq 1$. Then the second conditions imposes an (infinite, in
general) set of quadratic equations for these tensor
maps $\hat{\mu}_n$.

The natural isomorphism
$$
\odot^n(\fg[1])\simeq  (\Lambda^* \fg)[n]
$$
identifies the maps $\hat{\mu}_n: \odot^n (\fg[1]) \rar \fg[2]$
with the maps $\mu_n: \Lambda^n \fg \rar \fg[2-n]$,
$$
\hat{\mu}(v_1[1], \ldots, v_n[1]) = (-1)^{\sum_{i=1}^n
(n-i)\tv_i + n} \mu_n(v_1,\ldots v_n)[n].
$$
The homological condition $Q^2=0$ translates then precisely into the
higher Jacobi identities (\ref{id}) (see \cite{LS} for the proof).

\vst

%%%%%%%%%%%%%%%%%%%%%%%%%%%%%%%%%%%%%%%%
\paragraph{2.3. $L_{\infty}$-morphisms.} Given two
$L_{\infty}$-algebras, $(\fg, \mu_*)$ and $(\fg', \mu'_*)$. A
$L_{\infty}$-morphism $F$ from the first one to the second is, by
definition, a differential coalgebra homomorphism
$$
F: \left(\odot^*(\fg[1]), \Delta, Q\right) \lon \left(\odot^* (\fg'[1]), \Delta,
Q'\right).
$$
It is completely determined by a set of linear maps $\hat{F}_n: \odot^n (\fg[1])
\rar \fg'[1]$  (or, equivalently, by ${F}_n: \Lambda^n\fg \rar \fg'[1-n])$)
satisfying  an (infinite, in general) set of
equations. In the case when $(\fg', \mu'_n)$ is a differential
Abelian Lie algebra (i.e.\ $\mu'_n=0$ for $n\geq2$) these equations
take the form \cite{LM}
\Beq \label{abelian}
d F_n(v_1,\ldots, v_n)=\sum_{k+l=n+1}\sum_{\sigma\in
Sh(k,n)} (-1)^{\tilde{\sigma}+ k(l-1)} e(\sigma)
F_l\left(\mu_k(v_{\sigma(1)},\ldots,v_{\sigma(k)}),
v_{\sigma(k+l)},\ldots, v_{\sigma(n)}\right),
\Eeq
where we denoted $d=\mu_1'$.

\vst

A $L_{\infty}$-morphism $F:(\fg, \mu_*)\rar (\fg', \mu'_*)$ is
called a {\em quasi-isomorphism}\, if its first component $F_1: \fg\rar \fg'$
induces an isomorphism between cohomology groups of complexes
$(\fg, \mu_1)$ and $(\fg', \mu'_1)$. It is called a $L_{\infty}$-{\em
isomorphism}, if $F_1: \fg\rar \fh$ is an isomorphism of graded vector
spaces.

A differential Lie algebra $(\fg, d, [\, , \,])$ is
called {\em formal}\, if it is quasi-isomorphic to its cohomology
$(\rH^*(\fg), 0, [\, , \, ])$.
\vst

\noindent\mbox{{\bf 2.3.1. Fact} \cite{Ko3}}. If the
$L_{\infty}$-morphism $F:(\fg, \mu_*)\rar (\fg', \mu'_*)$ is
 a quasi-isomorphism, then there exists a $L_{\infty}$-morphism
 $F':(\fg', \mu'_*)\rar (\fg, \mu_*)$ which induces the inverse
 isomorphism between cohomology groups of complexes
$(\fg, \mu_1)$ and $(\fg', \mu'_1)$.

\vst

A similar statement holds true for $L_{\infty}$-isomorphisms.

\vst

%%%%%%%%%%%%%%%%%%%%%%%%%%%%%%%%%%%%%%%%%%%%%%%%%%%%
\paragraph{2.4. Maurer-Cartan equations.} Let $(\fg, \mu_*)$
be a $L_{\infty}$-algebra. A subset $\cM\cC(\fg)\subset \fg$
of {\em solutions to Maurer-Cartan equations}\, is defined formally as
follows\footnote{In our context $\cM\cC$ is applied only to
$L_{\infty}$-algebras of the form $\fg\otimes m$, where $m$
is the maximal ideal of an Artin algebra (or its completion);
hence no convergence problem arises.}

\begin{eqnarray*}
\cM\cC(\fg) &=& \left\{ \Gamma \in \fg^1 \mid
\sum_{k=1}^{\infty} \frac{(-1)^{k(k+1)/2}}{k!}
\mu_k(\Gamma, \ldots, \Gamma)=0\right\}\\
&=& \left\{ \Gamma \in \fg^1 \mid
d\Gamma +  \frac{1}{2!}\mu_2(\Gamma,\Gamma)
 - \frac{1}{3!}\mu_3(\Gamma,\Gamma,\Gamma)
- \frac{1}{4!}\mu_4(\Gamma,\Gamma,\Gamma,\Gamma) + \ldots =0
\right\}.
\end{eqnarray*}
If $\mu_n=0$ for $n\geq 3$, the equation on the r.h.s.\  reduces to the standard
Maurer-Cartan equation in a differential Lie algebra.

\vst

\noindent\mbox{{\bf 2.4.1. Fact} \cite{GM2,Ko3}}. Let
$F:(\fg, \mu_*)\rar (\fg', \mu'_*)$ be a $L_{\infty}$-morphism
between two $L_{\infty}$-algebras. If $F$ is such that $F_1$
provides an isomorphism of complexes
$(\fg, \mu_1)$ and $(\fg', \mu'_1)$, then $\cM\cC(\fg)=\cM\cC(\fg')$.
In particular, the deformation functors $\mbox{MDef}_{\fg}$ and
$\mbox{MDef}_{\fg,}$ are
equivalent. Moreover, the  last statement remains true under a weaker assumption that
 $F$ is  a quasi-isomorphism.

\vst

%%%%%%%%%%%%%%%%%%%%%%%%%%%%%%%%%%%%%%%%%%%%%%%%%%%%%%%
\paragraph{2.5. A geometric interpretation of a
$L_{\infty}$-algebra.} The dual of the free cocommutative coalgebra
$\odot^*(\fg[1])$
 can be identified with the algebra of
formal power series on the vector superspace $\fg[1]$ viewed as a
formal $\Z$-graded supermanifold (to emphasize this change of thought we denote
the supermanifold structure on $\fg[1]$ by $M_{\fg[1]}$). With
this identification the codifferential $Q$ on $\odot^*(\fg[1])$ goes
into a degree +1 vector field $Q$ on $M_{\fg[1]}$ satisfying the following two
conditions \cite{AKSZ,Ko3}
\begin{itemize}
\item[a)] $Q^2=0$;
\item[b)] $\left.Q\right|_0=0$.
\end{itemize}

The set $\cM\cC(\fg)$ is
then precisely the subset in $M_{\fg[1]}$ where $Q$ vanishes \cite{Ko3}. A
$L_{\infty}$-morphism $F$ between two $L_{\infty}$-algebras $(\fg,Q)$
and $(\fg',Q')$ is nothing but  a $Q$-equivariant map
between pointed formal graded supermanifolds $(M_{\fg[1]}, 0)$ and
$(M_{\fg'[1]},0)$. In this setting the Fact~2.4.1 becomes very
transparent: the morphism $F$ evidently maps zeros of $Q$ into
zeros of $Q'$.

\vst

Since $Q^2=0$, for any vector field $\alpha$ on $M_{\fg[1]}$ the
associated vector field $\left.[Q,\alpha]\right|_{\cM\cC(\fg)}$
is tangent to $\cM\cC(\fg)$. Moreover, since
$$
\left[[Q,\alpha], [Q,\beta]\right]=\left[Q, \left[\alpha,[Q,
\beta]\right]\right],
$$
such vector fields form an integrable distribution on
$\cM\cC(\fg)$ \cite{Ko3,Ba}. The leaves of this distribution define a natural
{\em gauge equivalence}\, on $\cM\cC(\fg)$.

\vst

We refer to the nice exposition of Kontsevich \cite{Ko3}
for more details on this geometrical model.

%%%%%%%%%%%%%%%%%%%%%%%%%%%%%%%%%%%%%%%%%%%%%%%%%%%%%%%%%%%%%%%%%%%
\paragraph{2.6. $\Gamma$-deformed $L_{\infty}$-structure.} Let $(\fg,
\mu_*)$,
or equivalently $(\fg, Q)$, be a $L_{\infty}$-structure on a
graded vector space $\fg$. There is an odd (more precisely, degree $+1$) linear morphism
$$
\Psi: \fg \lon \Gamma(M_{\fg[1]}, TM_{\fg[1]})
$$
identifying elements of $\fg$ with constant vector fields on the
supermanifold $M_{\fg[1]}$. In particular, an  (odd)
element $\Gamma$ of $\fg^1$
give rise to an (even)  constant vector field $\Psi(\Gamma)$
on $M_{\fg[1]}$ which in turn defines a  local diffeomorpism (the shift by $-\Gamma[1]$)
and hence
an associated morphism of Lie algebras
of formal vector fields,
$$
\begin{array}{rccc}
e^{ad \Psi(\Gamma)}:& \Gamma(M_{\fg[1]}, TM_{\fg[1]}) &\lon &
\Gamma(M_{\fg[1]}, TM_{\fg[1]})  \vspace{3mm} \\
& V & \lon & e^{ad \Psi(\Gamma)}(V)= V + \left[\Psi(\Gamma),V\right] +
\frac{1}{2!}\left[\Psi(\Gamma), [\Psi(\Gamma),V]\right] +
\ldots.
\end{array}
$$

\paragraph{2.6.1. Theorem.} {\em Let $(\fg, Q)$ be a
$L_{\infty}$-algebra. For any $\Gamma\in \cM\cC(\fg)$, the associated data
$$
(\fg,\,  Q_{\Gamma}:=e^{ad \Psi(\Gamma)}Q)
$$
is again a $L_{\infty}$-algebra}.
\vst

\noindent{\bf Proof}. The $\Gamma$-deformed vector field $Q_{\Gamma}$
is clearly homological,
$$
\left[Q_{\Gamma}, Q_{\Gamma}\right]=\left[e^{ad \Psi(\Gamma)}Q,
e^{ad \Psi(\Gamma)}Q\right] = e^{ad \Psi(\Gamma)}[Q,Q]=0.
$$
The main point is that the condition
$$
\left.Q_{\Gamma}\right|_0=0
$$
holds precisely when $\Gamma$ is a solution of the Maurer-Cartan
equations with respect to the original $L_{\infty}$-structure $Q$. $\Box$

\vst

The differential of the $\Gamma$-deformed $L_{\infty}$-structure
$(\fg, Q_{\Gamma})$
is given explicitly by
\begin{eqnarray}
d_{\Gamma}v &=& -\sum_{k=1}^{\infty} \frac{(-1)^{k(k+1)/2}}{(k-1)!}
\mu_k(\Gamma, \ldots, \Gamma, v)\nonumber\\
&=& d\Gamma +  \mu_2(\Gamma,v)
 - \frac{1}{2!}\mu_3(\Gamma,\Gamma,v)
- \frac{1}{3!}\mu_4(\Gamma,\Gamma,\Gamma,v) + \ldots,
\label{dgamma}
\end{eqnarray}
where $v\in \fg$.

%%%%%%%%%%%%%%%%%%%%%%%%%%%%%%%%%%%%%%%%%%%%%%%%%%%%%%%%%%%%%
\paragraph{2.7. $A_{\infty}$-algebras.}
A {\em strong homotopy algebra}, or shortly $A_{\infty}$-{\em algebra},
is by definition a $\Z$-graded vector space $A$ equipped with linear maps,
$$
\Ba{rcccc}
m_k: & \ot^k A & \lon & A[2-k] & \\
& a_1\ot\ldots \ot a_k & \lon & m_k(a_1, \ldots, a_k), & \ \ \ \ k\geq 1,
\Ea
$$
which satisfy, for any $n\geq 1$ and
 any $a_1, \ldots, a_n \in A$,
 the following
{\em higher order associativity conditions},
\Beq \label{idA}
\sum_{k+l=n+1}\sum_{j=0}^{ k-1} (-1)^r m_k\left(a_1,\ldots,a_j,
m_l(a_{j+1},\ldots, a_{j+l}), a_{j+l+1}, \ldots, a_n\right)=0,
\Eeq
where $r=\tilde{l}(\ta_1 +\ldots + \ta_j) + \tj(\tl-1) +
(\tk-1)\tl$.
\vst

It is easy to see from (\ref{idA}) that
$A_{\infty}$-algebras with $m_k=0$ for $k\geq 3$
are nothing but the usual graded differential associative algebras $(A,d, \cdot)$
with the differential $d=\mu_1$ and the associative multiplication given by
$a_1\cdot a_2=m_2(a_1,a_2)$.

\vst

There is a natural functor
$$
\Ba{ccc}
\left\{\Ba{l} \mbox{the category of $A_{\infty}$-structures}\\
                    \mbox{on a graded vector space $V$}
                    \Ea \right\}&
\lon & \left\{\Ba{l} \mbox{the category of $L_{\infty}$-structures}\\
                    \mbox{on a graded vector space $V$}\Ea
                    \right\}
\Ea
$$
given by \cite{LM}
$$
\Phi(m_n)(v_1, \ldots, v_n) = \sum_{\sigma in S_n}
(-1)^{\sigma}e(\sigma)
m_n(v_{\sigma(1)}, \ldots, v_{\sigma(n)}),
$$
where $S_n$ is the permutation group on $n$ elements. This is a
generalization of the usual construction of the Lie algebra out
of an associative algebra.

%%%%%%%%%%%%%%%%%%%%%%%%%%%%%%%%%%%%%%%%%%%%%%%%%%%%%%%%%%%%%%%
%%%%%%%%%%%%%%%%%%%%%%%%%%%%%%%%%%%%%%%%%%%%%%%%%%%%%%%%%%%%%%%
\section{An unobstructed deformation functor}
%%%%%%%%%%%%%%%%%%%%%%%%%%%%%%%%%%%%%%%%%%%%%%%%%%%%%%%%%%%%%%
\paragraph{3.1. Theorem.}  {\em Let $(\fg, d, [\, , \,])$ be a
differential graded Lie algebra, and let $\eta: \fg\rar \fg[-1]$ be any linear map.
Then the formulae
\begin{eqnarray}
\mu_1(v_1)&:=& dv_1 \nonumber\\
\mu_2(v_1,v_2) &:=& (d\eta+\eta d)[v_1,v_2] \nonumber\\
\mu_3(v_1,v_2,v_3) &:=& -\eta\left[\mu_2(v_1,v_2), v_3\right] +
(-1)^{\tv_2\tv_3}\eta\left[\mu_2(v_1,v_3), v_2\right]
- (-1)^{\tv_1(\tv_2+\tv_3)}\eta\left[\mu_2(v_2,v_3), v_1\right] \nonumber\\
\cdots && \nonumber\\
\mu_n(v_1,\ldots,v_n)&:=& (-1)^n \sum_{\sigma\in Sh(n-1,n)}
(-1)^{\tilde{\sigma}}e(\sigma) \eta\left[\mu_{n-1}(v_{\sigma(1)}, \ldots,
v_{\sigma(n-1)}), v_{\sigma(n)}\right] \label{linf}\\
\cdots && \nonumber
\end{eqnarray}
define inductively a structure of $L_{\infty}$-algebra on the graded vector
space $\fg$. }

\vst

\noindent{\bf Proof}~(an outline). We have to show that the
tensors
\Beqrn
\Phi_{n}(v_1,\ldots,v_n)&:=&\sum_{k+l=n+1}\sum_{\sigma\in Sh(k,n)}
(-1)^{\tilde{\sigma}+k(l-1)} e(\sigma)
 \mu_l\left(\mu_k(v_{\sigma(1)},\ldots,v_{\sigma(k)}),
v_{\sigma(k+l)},\ldots, v_{\sigma(n)}\right)\\
&=& d\mu_n(v_1, \ldots, v_n) + \sum_{i=1}^n (-1)^{n-1 + \tv_1 + \ldots +
\tv_{i-1}}\mu_n(v_1, \ldots, v_i, dv_i, v_{i+1}, \ldots, v_n) \\
&& + (-1)^{n-1} \sum_{\sigma\in Sh(n-1,n)}(-1)^{\tilde{\sigma}}
e(\sigma)
\mu_2\left(\mu_{n-1}(v_{\sigma(1)},\ldots,v_{\sigma(n-1)}),
v_{\sigma(n)}\right)\\
&& + \sum_{k+l=n+1 \atop k+1,l\geq 3}\sum_{\sigma\in Sh(k,n)}
(-1)^{\tilde{\sigma}+k(l-1)} e(\sigma)
 \mu_l\left(\mu_k(v_{\sigma(1)},\ldots,v_{\sigma(k)}),
v_{\sigma(k+l)},\ldots, v_{\sigma(n)}\right)
\Eeqrn
constructed out of the maps $\mu_*$ defined above vanish for all $n\geq
2$ and all $v_1, \ldots,v_n\in \fg$.

\vst

For $n\geq 3$,
\Beqrn
d\mu_n(v_1,\ldots,v_n)&=& (-1)^n \sum_{\sigma\in Sh(n-1,n)}
(-1)^{\sigma}e(\sigma)d \eta\left[\mu_{n-1}(v_{\sigma(1)}, \ldots,
v_{\sigma(n-1)}), v_{\sigma(n)}\right]\\
 &=& (-1)^n \sum_{\sigma\in Sh(n-1,n)}
(-1)^{\sigma}e(\sigma) \mu_2\left[\mu_{n-1}(v_{\sigma(1)}, \ldots,
v_{\sigma(n-1)}), v_{\sigma(n)}\right]\\
&& + (-1)^{n+1} \sum_{\sigma\in Sh(n-1,n)}
(-1)^{\sigma}e(\sigma) \eta d\left[\mu_{n-1}(v_{\sigma(1)}, \ldots,
v_{\sigma(n-1)}), v_{\sigma(n)}\right].
\Eeqrn
Substituting this expression into $\Phi_{n}(v_1, \ldots,v_n)$, one
gets,
after tedious but straightforward algebraic manipulations, a recursion
formula
$$
\Phi_n(v_1, \ldots, v_n) = (-1)^{n+1} \sum_{\sigma\in Sh(n-1,n)}
(-1)^{\tilde{\sigma}}e(\sigma) Q\left[ \Phi_{n-1}(v_{\sigma(1)}, \ldots,
v_{\sigma(n-1)}), v_{\sigma(n)}\right], \ \ \ n\geq 3.
$$

Since $\Phi_1(v_1)=d^2v_1=0$ and
\Beqrn
\Phi_2(v_1,v_2) &=& d\mu_2(v_1,v_2) - \mu_2(dv_1,v_2) -
(-1)^{\tv_1}\mu_2(v_1,dv_2)\\
&=& d\eta d [v_1,v_2] - d\eta [dv_1, v_2]
- (-1)^{\tv_1}d\eta[v_1,dv_2]\\
&=& 0,
\Eeqrn
the required statement follows. $\Box$

\paragraph{3.1.1. Remark.} Let $(\fg,d,[\, , \, ], \eta:\fg\rar \fg )$ be
the same data as in Theorem~3.1. Setting formally $\lambda_1 := -
\eta^{-1}$ we define a series of linear maps,
$$
\la_n: \Lambda^n \fg \lon \fg[2-n], \ \ \ \ \ n\geq 2,
$$
by a recursive formula (cf.\ \cite{Me1})
$$
\la_n(v_1,\ldots,v_n) :=   \sum_{ k+l=n\atop k,l\geq 1
}\sum_{\sigma\in Sh(k,n)}(-1)^{\tilde{\sigma}+r}e(\sigma)
\left[\eta\la_k(v_{\sigma(1)}, \ldots, v_{\sigma(k)}),
\eta\la_l(v_{\sigma(k+1)},\ldots,v_{\sigma(n)})\right],
$$
where $r={k+1 +(l-1)(\tv_{\sigma(1)}+\ldots+ \tv_{\sigma(k)})}$.
Then the data
\Beqrn
m_1 &:=& d\\
m_n &:=& (1-[d,\eta])\la_n, \ \ \ \mbox{for}\ n\geq 2,
\Eeqrn
define a structure of $L_{\infty}$-algebra on $\fg$
which is, in a sense, complementary to the one given in Theorem~3.1.
We do not use this structure in the paper and hence omit the proof.

\paragraph{3.1.2. Remark.} The Theorem~3.1 has an $A_{\infty}$
analogue.
Let $(A,d,\cdot)$  be a differential graded associative algebra and
let $\eta:A\rar A[-1]$ be any linear operator. Then the formulae
\begin{eqnarray*}
m_1(a_1)&:=& dv_1\\
m_2(a_1,a_2) &:=& (d\eta+\eta d)(a_1\cdot a_2)\\
m_3(a_1,a_2,a_3) &:=& \eta\left(-m_2(a_1,a_2)\cdot a_3 +
a_1\cdot m_2(a_2,a_3)\right)\\
\cdots &&\\
m_n(a_1,\ldots,a_n)&:=& \eta\left((-1)^n
m_{n-1}(a_1, \ldots, a_{n-1})\cdot a_{n}
+(-1)^{(n-1)\ta_1} a_1 \cdot m_{n-1}(a_2, \ldots, a_{n})
\right)\\
\cdots &&
\end{eqnarray*}
define a structure of $A_{\infty}$-algebra on the vector
superspace $A$. We omit the proof.

%%%%%%%%%%%%%%%%%%%%%%%%%%%%%%%%%%%%%%%%%%%%%%%%%%%%%%%%%%%%%%%%%%%%%
\paragraph{3.2. Theorem.} {\em
Let $(\fg, d, [\, , \,])$ be a
differential graded Lie algebra and  $\eta: \fg\rar \fg[-1]$  any  linear
map.
Then the {\em Kuranishi map}\, $K:\Lambda^* \fg \rar \fg[1-*]$
given by its homogeneous components as follows
\Beqrn
K_1(v_1) &:=& v_1, \\
K_2(v_1,v_2)&:=& \eta[v_1,v_2],\\
K_n(v_1,\ldots,v_n)&:=& 0 \hspace{2cm} \mbox{for}\ n\geq 3,
\Eeqrn
defines a $L_{\infty}$-isomorphism between
the $L_{\infty}$-structure $(\fg, \mu_*)$ induced on $\fg$ by
Theorem~3.1
and the differential Abelian Lie algebra $(\fg, d, 0)$.}

\vst

\noindent{\bf Proof.} The equations (\ref{abelian}) defining
$L_{\infty}$-morphisms into an a differential Abelian Lie algebra
take, in our case, the form
$$
\Ba{ll}
n=1: & dv_1 - dv_1=0,\\
n=2: & dQ[v_1,v_2] - \mu_2(v_1,v_2) + Q[dv_1,v_2] + (-1)^{\tv_1}Q[v_1,
dv_2] =0, \\
n\geq 3: & \mu_n(v_1, \ldots, v_n) -
(-1)^n \sum_{\sigma\in Sh(n-1,n)}
(-1)^{\tilde{\sigma}}e(\sigma) \eta\left[\mu_{n-1}(v_{\sigma(1)}, \ldots,
v_{\sigma(n-1)}), v_{\sigma(n)}\right] =0.
\Ea
$$
They are all obviously satisfied. $\Box$

\vst

Therefore, the ``naive'' $L_{\infty}$-structure (\ref{linf}) on $\fg$ is precisely the one
obtained from the differential Abelian Lie algebra $(\fg,d,0)$ by
inverting (in the category of $L_{\infty}$-algebras) the Kuranishi
\cite{Ku1,Ku2} map. This is a key observation of the paper from
which everything else follows.

%%%%%%%%%%%%%%%%%%%%%%%%%%%%%%%%%%%%%%%%%%%%%%%%%%%%%%%%%%%%%%%
\paragraph{3.3. Corollary.}
It follows from Theorem~3.2 and
Fact~2.4.1 that the deformation functor
$$
\Ba{rccc}
\mbox{MDef}^{\Z}_{\fg}: & \left\{\Ba{l} \mbox{the category of}\\
                    \mbox{graded Artin}\\
                    \mbox{$k$-local algebras}\Ea \right\}&
\lon & \left\{\mbox{the category of sets}\right\} \\
& \cB & \lon & \frac{
\left\{\Gamma \in (\fg\ot
m_{\cB})^{{1}} \mid \sum_{k=1}^{\infty} \frac{(-1)^{k(k+1)/2}}{k!}
\mu_k(\Gamma, \ldots, \Gamma)=0\right\}}{\mbox{natural gauge equivalence}}
\Ea
$$
is equivalent to the deformation functor
$$
\mbox{Def}_{\fg}: \cB  \lon  \frac{
\left\{\Gamma \in (\fg\ot
m_{\cB})^{{1}} \mid d\Gamma =0\right\}}{\mbox{Im}\, d},
$$
which is pro-representable by the graded algebra $k[[t_{\bf H^*}]]$
of formal power series on the $[1]$-shifted cohomology space,
$$
{\bf H^*}\equiv \rH^*(\fg,d)[1]:=\frac{\mbox{Ker}\, d}{\mbox{Im}\, d}[1],
$$
viewed as a graded linear supermanifold. Hence $\mbox{MDef}^{\Z}_{\fg}$
is always pro-representable by a germ of a  pointed formal supermanifold
$(\cM^{\Z}_{\fg}, p)$ which is isomorphic to the completion of the
analytic germ $\f_0$ of the pointed supermanifold $({\bf H^*}, 0)$.

\vst

Similarly, the deformation functor $\mbox{MDef}^{0}_{\fg}$ (defined in
the Subsect.\ 1.3)
is always pro-representable by the algebra
of formal power series on the vector superspace ${\bf H^1}\equiv
\Pi\rH^1(\fg,d)$, $\Pi$ being the parity change functor.
The associated formal pointed manifold we denote by
$(\cM_{\fg}^0, p)$.

%%%%%%%%%%%%%%%%%%%%%%%%%%%%%%%%%%%%%%%%%%%%%%%%%%%%%%%%%%%%%%%%%%%%
\paragraph{3.4. A versal solution to the Maurer-Cartan equations.}
Let us fix a basis $[\gamma_{\al}]$ of $\rH^*(\fg,d)$. It defines an associated
basis $\ga_{\al}[1]$ in ${\bf H^*}$ and hence a set of linear coordinates $\{t^{\al}\}$.
Let $\ga_{\al}\in \mbox{Ker}\, d$
be any representatives of the cohomology classes $[\ga_{\al}]$ in $\fg$.
Then
$$
\tilde{\Gamma}(t):= \sum_{\al} \ga_{\al} t^{\al}
$$
is a versal solution to the Maurer-Cartan equation in the Abelian Lie
algebra $(\fg\ot k[[t_{\bf H^*}]],d,0)$. By Theorem~3.3 and Fact 2.4.1
any formal power series
$$
\Gamma(t)= \sum_{\al} t^{\al} \Gamma_{\al} + \sum_{\al_1, \al_2}\Gamma_{\al_1,{\al_2}}
t^{\al_1} t^{\al_2} + \ldots \in (\fg \ot k[[t_{\bf H^*}]])^1
$$
satisfying the equation
\Beqr
\tilde{\Gamma}(t) &=& \sum_{n=1}^{\infty} \frac{1}{n!} K_n(\Gamma(t), \ldots,
\Gamma(t)) \nonumber\\
&=& \Gamma(t) + \frac{1}{2} \eta [\Gamma, \Gamma] \label{kur}
\Eeqr
gives a versal solution to the Maurer-Cartan equation,
\Beq
\sum_{k=1}^{\infty} \frac{(-1)^{k(k+1)/2}}{k!}
\mu_k(\Gamma(t), \ldots, \Gamma(t))=0,  \label{MC}
\Eeq
in the $L_{\infty}$-algebra $(\fg\ot k[[t_{\bf H^*}]], \mu_*)$.
The equation (\ref{kur}) is easily solved by
$$
\Gamma(t)= \Gamma_1(t) + \Gamma_2(t) + \ldots + \Gamma_n(t) + \ldots
$$
where
\Beqr
\Gamma_1 &=& \sum_{\al} \ga_{\al} t^{\al}\nonumber \\
\Gamma_2 &=& - \frac{1}{2} \eta [\Gamma_1(t), \Gamma_1(t)],\nonumber\\
\Gamma_3 &=&  - \frac{1}{2} \eta \left([\Gamma_1(t), \Gamma_2(t)] +
[\Gamma_2(t), \Gamma_1(t)]\right),\nonumber\\
\ldots && \nonumber\\
\Gamma_n &=& - \frac{1}{2} \eta \left(\sum_{k=1}^{n-1}[\Gamma_k(t),
\Gamma_{n-k}(t)]
\right) \label{versal}\\
\ldots && \nonumber
\Eeqr
This is a well known  power series  \cite{Ko,Ku1,Ku2, Kob}  playing a key role
in the deformation theory.  Therefore, it has a very simple algebraic
interpretation  within the category of $L_{\infty}$-algebras.

\vst

\paragraph{3.4.1. Remark.}
The same formulae (\ref{versal})
describe a versal solution of the Maurer-Cartan
equations in the $L_{\infty}$-algebra $(\fg\ot k[[t_{\bf H^1}]], \mu_*)$.

%%%%%%%%%%%%%%%%%%%%%%%%%%%%%%%%%%%%%%%%%%%%%%%%%%%%%%%%%%%%%%%%%%%%%%%%%%
\paragraph{3.5. Hodge structures.} Assume that the graded differential Lie
algebra $(\fg, d, [\, , \, ])$ is equipped with a norm $|| \ ||_i$ on
each $\fg^i\subset \fg$ making $\fg^i$ into a normed vector space such
that both the differential $d: \fg^i \rar \fg^{i+1}$ and the Lie bracket
$[\, , \, ]:\fg^i \times \fg^j\rar \fg^{i+j}$ are continuous. Assume also that
there exists a  projection
$$
P_{\bf H}: \fg \lon {\bf H}
$$
and a linear operator $\eta:\fg \lon \fg[-1]$ such that
the Hodge decomposition holds
$$
\mbox{Id}= P_{\bf H} + d\eta  + \eta d,
$$
and similarly for the completion of $\fg$
with respect to the norm. In the usual Hodge theory $\eta=d^* G$,
where $d^*$ is the adjoint to $d$ and $G$ is the Green function.

\paragraph{3.5.1. Smoothness.}
Using the implicit function theorem in a Banach space
as in \cite{GM2}, one easily shows that
 the pointed formal supermanifold $({\cM}^{\Z}_{\fg},p)$
representing the functor $\mbox{MDef}^{\Z}_{\fg}$ has a
natural smooth analytic structure which makes $({\cM}^{\Z}_{\fg},p)$  analytically
diffeomorphic (with respect to the Kuranishi map) to a neighbourhood of zero
in  the vector space ${\bf H^*}$ (it is also easy to show
that, for sufficiently small $t$, the power series
(\ref{versal}) is convergent \cite{Ko,Kob}).
This proves the Claim (iii)
in the Introduction.

\vst

In a similar  way one provides the pointed moduli space
$({\cM}^{0}_{\fg}, p)$ representing the functor $\mbox{MDef}^{0}_{\fg}$ with a smooth
analytic structure in such a way that the Kuranishi map becomes an
analytic equivalence between $({\cM}^{0}_{\fg}, p)$ and $({\bf H^1},
0)$. Hence the Claim~(i) follows.

%%%%%%%%%%%%%%%%%%%%%%%%%%%%%%%%
\paragraph{3.5.2. Vanishing obstructions.}
Assume that in a formal differential Lie
algebra $(\fg, d, [\, , \, ])$ the induced map
\Beq
[\, , \, ]: \rH^*(\fg) \times \rH^*(\fg) \lon \rH^*(\fg), \label{obs2}
\Eeq
is zero. Put another way, for any $v_1$ and $v_2$ in $\fg$,
$$
P_{\bf H}[P_{\bf H}(v_1), P_{\bf H}(v_2)] =0.
$$
Then there is a chain of quasi-isomorphisms
$$
(\fg, \mu_*) \stackrel{{\small Kuranishi}\atop{{\small map}}}{\lon} (\fg,d,0) \stackrel{P_{\bf H}}{\lon} ({\bf H},
0,0)
{\longleftarrow} (\fg, d, [\, , \, ]),
$$
connecting the $L_{\infty}$-structure (\ref{linf}) with the original
differentiable Lie algebra structure on $\fg$. By Fact~2.4.1,
the deformation functors $\mbox{MDef}_{\fg}^{\Z}$ and
$\mbox{Def}_{\fg}^{\Z}$ are equivalent. This  in turn implies
that the completion of the  analytic germ of the pointed
moduli space $({\cM}^{\Z}_{\fg},p)$ (constructed in Subsect.\ 3.5.1)
is isomorphic to the completion
  of the analytic germ of the classical Kuranishi pointed moduli space
$(\cK^{\Z}_{\fg},p)$ (constructed in \cite{Ku1,Ku2,GM2}).
By \cite{Artin}, the analytic
equivalence of $({\cM}^{\Z}_{\fg},p)$ and $({\cK}^{\Z}_{\fg},p)$ follows.
 This completes the proof of the
Claim (iv) in the Introduction.

\paragraph{3.5.3. Non-vanishing obstructions.}  The classical Kuranishi
moduli space $\cK^0_{\fg}$ is an analytic subspace in ${\cM}^0_{\fg}$
given by the equations \cite{Ku1,Ku2}
$$
P_{\bf H}[\Gamma(t), \Gamma(t)]=0,
$$
where $\Gamma(t)\in \fg^1 \ot k[[t_{\bf H^1}]]$ is the versal solution
(\ref{versal}) (see  Remark~3.4.1).

\vst

For any $t\in \cK^0_{\fg}$ one has
\Beqrn
\mu_2(\Gamma(t), \Gamma(t)) &=& (d\eta + \eta d)[\Gamma(t), \Gamma(t)]\\
& =& (1- P_{\bf H})[\Gamma(t), \Gamma(t)] \\
&=& [\Gamma(t), \Gamma(t)],\\
\mu_3(\Gamma(t), \Gamma(t),\Gamma(t)) &=&
 -3\eta\left[[\Gamma(t), \Gamma(t)], \Gamma(t)\right] =0,\\
\Eeqrn
and hence
$$
\mu_n(\Gamma(t),\ldots,\Gamma(t))= (-1)^n n
\eta\left[\mu_{n-1}(\Gamma(t), \ldots, \Gamma(t)), \Gamma(t)\right]=0
$$
for all $n\geq 3$. Therefore, the classical Kuranishi space $\cK^0_{\fg}$
is precisely the subspace of  ${\cM}^0_{\fg}$ where the Maurer-Cartan
equation (\ref{MC}) of the $L_{\infty}$-algebra (\ref{linf})
degenerates into the following one,
$$
d\Gamma(t) + \frac{1}{2}[\Gamma(t), \Gamma(t)] =0.
$$
This is exactly the Maurer-Cartan equation of the differential Lie
algebra $(\fg, d, [\, , \, ])$. This explains  the Claim (ii) in the
Subsect.\ 1.3.

%\vst

\paragraph{\bf 3.6. Open questions and speculations.}
If $\rH^2(\fg,d)=0$, then, evidently, ${\cM}^0_{\fg}\simeq {\cK}^0_{\fg}$
and the deformed differential (\ref{dgamma})
degenerates into a usual linear connection
$$
d_{\Gamma} = d + \mbox{ad} \Gamma.
$$
If obstructions do not vanish, then $d_{\Gamma}$ has terms higher
order in $\Gamma$. Though it is still ``flat", $d_{\Gamma}^2=0$, this
object is no more a linear connection. It seems that, for a
geometric interpretation of the smooth moduli space $\cM^{\Z}_{\fg}$
in the case
$(\fg, d,  [\, , \,])=(E\ot E^* \ot \Omega^{0,\bullet}M, \bar{\partial}, \mbox{standard
bracket})$, $E$ being a holomorphic vector bundle over a complex
manifold $M$, one should switch from the category of projective modules
 to the category of strong homotopy modules
over the differential algebra
$(\Omega^{0,\bullet}M, \bar{\partial})$ (cf.\ \cite{Kapr}) or
even over its $A_{\infty}$-versions.

\vst

If $M$ is a Calabi-Yau manifold, then the moduli space $\cM^{\Z}_{\fg}$
associated with the differential graded algebras $(\Lambda^{\bullet}TM\ot
\Omega^{0,\bullet}M, \bar{\partial}, \mbox{Schouten bracket})$
is precisely the Barannikov-Kontsevich \cite{BaKo} extended moduli space
of complex structures. To understand $\cM^{\Z}_{\fg}$ for a general
compact
complex manifold, one should probably think about a strong homotopy
generalization of the notion of odd contact structure on a complex
supermanifold.

\vse

%\pagebreak
{\small

{\small
\begin{tabular}{l}
Department of Mathematics\\
University of Glasgow\\
15 University Gardens \\
 Glasgow G12 8QW, UK\\
 \mbox{\rm sm@maths.gla.ac.uk}
\end{tabular}
}

\end{document}